\documentclass[12pt]{amsart}
\usepackage[centertags]{amsmath}
\usepackage{amsfonts,amssymb}
\usepackage{graphicx}
\usepackage{enumerate}
\usepackage[latin1]{inputenc}
\usepackage{float}

  \newtheorem{theorem}{Theorem}
  \newtheorem{corollary}{Corollary}
  \newtheorem{proposition}{Proposition}
  \newtheorem{lemma}{Lemma}%
  \theoremstyle{remark}
  
\newtheorem{conjecture}{Conjecture}

 \newcommand{\perm}{{\rm perm}}
 \newcommand{\oo}{{\rm o}}

\begin{document}

\title{Coprime permutations}
\date{\today}

\author{Carl Pomerance}
\address{Mathematics Department, Dartmouth College, Hanover, NH 03784}
\email{carlp@math.dartmouth.edu}

\begin{abstract}
Let $C(n)$ denote the number of permutations $\sigma$ of $[n]=\{1,2,\dots,n\}$
such that $\gcd(j,\sigma(j))=1$ for each $j\in[n]$.  We prove that for $n$
sufficiently large, $n!/3.73^n < C(n) < n!/2.5^n$.

\end{abstract}

\subjclass[2010]{11A25, 11B75, 11N60}

\keywords{coprime permutation, coprime matching, distribution function, Euler's function}
\maketitle

\vskip-30pt
\newenvironment{dedication}
        {\vspace{6ex}\begin{quotation}\begin{center}\begin{em}}
        {\par\end{em}\end{center}\end{quotation}}
\begin{dedication}
{In memory of Andrzej Schinzel (1937--2021)}
\end{dedication}
\vskip20pt

\section{Introduction}

Several papers, some recent, have dealt with coprime matchings between
two sets of $n$ consecutive integers; that is a matching where corresponding
pairs are coprime.  For example in a paper \cite{PS} with Selfridge, we showed
such a matching always exist if one of the intervals is $[n]=\{1,2,\dots,n\}$.
In Bohman and Peng \cite{BP} it is shown that a matching always exists if
$n$ is even and the numbers involved are not too large as a function of $n$,
with an interesting application to the lonely runner problem in Diophantine
approximations.  Their result was somewhat strengthened in \cite{P}.

The current paper considers the situation when both intervals are $[n]$.
In this case it is trivial that a coprime matching exists, just take the cyclic
permutation $(1,2,\dots,n)$.  So instead we consider the enumeration
problem.  Let $C(n)$ denote the number of permutations $\sigma$ of $[n]$
where $\gcd(j,\sigma(j))=1$ for each $j\in[n]$.  This problem was considered
in Jackson \cite{J} where $C(n)$ was enumerated for $n\le24$.  For
example,
\[
C(24)=1{,}142{,}807{,}773{,}593{,}600.
\]
After factoring his values, Jackson notes the appearance of sporadically large
primes, which indicates there may not be a simple formula.  The sequence
also has an OEIS page, see \cite{O}, where the value of $C(25)$, due to
A. P. Heinz, is presented (and the value for $C(16)$ is corrected).  There
are also links to further computations, especially those of Locke.  In Section
\ref{sec:comp} we discuss how $C(n)$ can be computed and verify Locke's
values.

Our principal result is the following.
\begin{theorem}
\label{thm:main}
For all large $n$, $n!/3.73^n < C(n) < n!/2.5^n$.
\end{theorem}

Important in the proof of the lower bound is a numerically explicit estimation
of the distribution function for $\varphi(n)/n$, where $\varphi$ is Euler's function.

It would seem likely that there is a constant $c$ with $2.5\le c\le 3.73$ with
$C(n)=n!/(c+o(1))^n$ as $n\to\infty$.  In Section \ref{sec:ub} we give some
thoughts towards this possibility.

In Section \ref{sec:comp} we discuss the numerical calculation of $C(n)$.
Finally, in Section \ref{sec:anti} we briefly discuss the number of permutations
$\sigma$ of $[n]$ where $\sigma(1)=1$ and for $2\le j\le n$, $\gcd(j,\sigma(j))>1$.

\section{Preliminaries}
\label{sec:pre}

Regarding notation, we have
\[
[n]=\{1,2,\dots,n\},\quad [n]_\oo=\{1,3,\dots,2n-1\}.
\]
Thus, $[n]_\oo$ is the set of the first $n$ odd positive integers.
Let $C_0(n)$ denote the number of one-to-one functions 
\[
f: [n]_\oo\longrightarrow[n]
\]
 such that each $\gcd(i,f(i))=1$.  Similarly, let
$C_1(n)$ denote the number of one-to-one functions 
\[
f:[n]\longrightarrow [n+1]_\oo
\]
 such that each $\gcd(i,f(i))=1$.
\begin{lemma}
\label{lem:odd}
We have $C(2n)=C_0(n)^2$ and for $n\ge2$,
$2C_0(n-1)^2\le C(2n+1)\le C_1(n)^2$.
\end{lemma}
\begin{proof}
Let $\sigma$ be a coprime permutation of $[2n]$.  Then $\sigma$ maps evens
to odds and odds to evens, so that $\sigma$ corresponds to a pair of coprime
matchings $\sigma_0,\sigma_1$ where $\sigma_0$ maps $\{2,4,\dots,2n\}$
to $\{1,3,\dots,2n-1\}$ and $\sigma_1$ maps $\{1,3,\dots,2n-1\}$ to
$\{2,4,\dots,2n\}$.  Then $f(2i-1)=\frac12\sigma_1(2i-1)$ is one of the
maps counted by $C_0(n)$ and so is $g(2i-1)=\frac12\sigma_0^{-1}(2i-1)$.
Conversely, each such pair of maps
corresponds to a coprime permutation $\sigma$ of $[2n]$.  This proves
that $C(2n)=C_0(n)^2$.

The upper bound for $C(2n+1)$ follows in the same way.  Let $\sigma$
be a coprime permutation of $[2n+1]$ and let $\sigma_0$ be $\sigma$ restricted
to even numbers. 
Then define $f_\sigma(i)=\sigma_0(2i)$, so that $f_\sigma$ is one of the functions counted
by $C_1(n)$.  Note that there is some $a\in\{1,3,\dots,2n+1\}$ with
$\sigma(a)$ odd, but all other members $b$ of $\{1,3,\dots,2n+1\}$ have
$\sigma(b)$ even.  Let $\sigma_1$ be $\sigma$ restricted to 
$S_{a,n}:=[n+1]_\oo\setminus\{a\}$ and let $g_\sigma(2i-1)=\frac12\sigma_1(2i-1)$
for $i\in S_{a,n}$.  Then $g_\sigma^{-1}$ is one of the functions counted by
$C_1(n)$.   Note that if $\tau$ is a coprime permutation of $[2n+1]$ such that
$f_\tau=f_\sigma$ and $g_\tau=g_\sigma$, then $\tau=\sigma$.  
This proves that $C(2n+1)\le C_1(n)^2$.  Note that the proof
ignores the condition $\gcd(a,\sigma(a))=1$, so it only gives an upper bound.

For the lower bound, note that $C(2n+1)\ge 2C(2n-2)$.  Indeed, 
corresponding to a coprime
permutation of $[2n-2]$ we augment it with either the cycle
$(2n-1,2n,2n+1)$ or its inverse, giving two coprime permutations of $[2n+1]$.
The lower bound in the lemma for $C(2n+1)$ now follows from the first part
of the lemma.
\end{proof}

 We remark that the sequence $C(1),C(2),\dots$ is not monotone, but it is
 monotone restricted to integers of the same parity.  Indeed, augmenting
 a coprime permutation of $[n]$ with the cycle $(n+1,n+2)$ gives a coprime
 permutation of $[n+2]$, so that $C(n)\le C(n+2)$.

Since $C_0(n)\le n!$ and $C_1(n)\le (n+1)!$, Lemma \ref{lem:odd} immediately
gives us that $C(2n)\le (n!)^2$ and $C(2n+1)\le (n+1)!^2$.  With Stirling's
formula this gives $C(n)\le n!/(2+o(1))^n$ as $n\to\infty$.  Note that this
argument considers only parity.  By bringing in 3, 5, etc., we can improve
this upper bound.  In Section \ref{sec:ub} we begin this process and show that
$C(n) <n!/(5/2)^n$ for all large $n$.

It is much harder to get a comparable lower bound for $C(n)$, and this is
our undertaking in the next two sections.  From the thoughts above it suffices
to get a lower bound for $C_0(n)$.
The lower bound in Theorem \ref{thm:main} is a consequence of the following result.

\begin{theorem}
\label{thm:main2}
For all large $n$, 
$C_0(n)\ge n!/1.864^n$.
\end{theorem}

\section{The distribution function}

Let $\omega(n)$ denote the number of distinct prime factors of $n$.

\begin{lemma}
\label{lem:sieve}
For positive integers $m,n$, the number
of  $j\le n$ with $\gcd(j,m)=1$ is within $2^{\omega(m)-1}$ of  $(\varphi(m)/m)n$.
\end{lemma}
\begin{proof}
The result is clear if $m=1$, so assume that $m>1$.
With $\mu$ the M\"obius function, the exact number of $j$'s is
\[
\sum_{\substack{d|m}}\sum_{\substack{j\le n\\d|j}}\mu(d)
=\sum_{\substack{d|m}}\left(\frac{\mu(d)}dn+\mu(d)\theta_d\right),
\]
where $0\le\theta_d<1$.  The sum of the main terms is $(\varphi(m)/m)n$.
There are $2^{\omega(m)}$ error terms $\mu(d)\theta_d$ with $\mu(d)\ne0$,
and since $m>1$, half of them are $\ge0$ and half are $\le0$. 
So the sum of the error terms has absolute magnitude $<2^{\omega(m)-1}$.
\end{proof}
\begin{corollary}
\label{cor:sieve}
For $m\le n$, the number of $j\le n$ with $\gcd(j,m)=1$ is greater than
$(\varphi(m)/m)n-\sqrt{n}$.
\end{corollary}
\begin{proof}
A short induction argument shows that
$\sqrt{m}>2^{\omega(m)-1}$, so the result follows directly from the lemma.
\end{proof}

\begin{lemma}
\label{lem:2ndmom}
For all large $n$ we have
\[
\sum_{\substack{m<2n\\m\,{\rm odd}}}\left(\frac m{\varphi(m)}\right)^2<1.78n.
\]
\end{lemma}
\begin{proof}
Define a multiplicative function $h$ with $h(p)=(2p-1)/(p-1)^2$ 
for each prime $p$ and $h(p^a)=0$ for $a\ge2$.  Then
\begin{align*}
\sum_{\substack{m<2n\\m\,{\rm odd}}}\left(\frac m{\varphi(m)}\right)^2
&=\sum_{\substack{m<2n\\m\,{\rm odd}}}\sum_{d\mid m}h(d)\\
&=\sum_{\substack{d<2n\\d\,{\rm odd}}}h(d)\sum_{\substack{j<2n/d\\j\,{\rm odd}}}1
=\sum_{\substack{d<2n\\d\,{\rm odd}}}h(d)\left(\frac nd+O(1)\right).
\end{align*}
The main term is
\[
<n\prod_{p>2}\left(1+\frac{2p-1}{(p-1)^2p}\right),
\]
and this infinite product converges to a constant smaller than $1.7725$.
For the error term a simple calculation shows that it is $O((\log n)^2)$,
so our conclusion follows.
\end{proof}

Let $\delta_\varphi(\alpha)$ be the distribution function for $\varphi(m)/m$;
that is, for $0\le\alpha\le1$,
\[
\delta_\varphi(\alpha)=\lim_{n\to\infty}\frac1n\sum_{\substack{m\le n\\\varphi(m)/m\le\alpha}}1.
\]
It is known after various papers of Schoenberg, Behrend, Chowla, Erd\H os, and Erd\H os--Wintner
 that the limit exists, $\delta_\varphi(0)=0$, $\delta_\varphi(1)=1$, and $\delta_\varphi$ is
strictly increasing and continuous.  In addition, at a dense set of numbers
in $[0,1]$, namely the values of $\varphi(m)/m$, the distribution function $\delta_\varphi$
has an infinite left
derivative.  This all can be generalized to odd numbers.
For $0\le\alpha\le 1$, let $D(\alpha,n)$ denote the number of odd $m<2n$
with $\varphi(m)/m\le \alpha$.  As with $\delta_\varphi$,
\[
\delta(\alpha):=\lim_{n\to\infty}D(\alpha,n)/n
\]
 exists, with $\delta$
continuous and strictly increasing on $[0,1]$, with $\delta(0)=0$ and $\delta(1)=1$.
By extending it to take the value 1 when $\alpha>1$, we have
\[
\delta_\varphi(\alpha)=\frac12(\delta(\alpha)+\delta(2\alpha)),
\]
as noted in \cite{PS}.  In particular, for $\frac12\le\alpha\le1$,
\begin{equation}
\label{eq:dist}
\delta(\alpha)=2\delta_\varphi(\alpha)-1.
\end{equation}

A consequence of the argument in \cite{PS} is that $\delta(\alpha)\le\alpha$
on $[0,1]$.  We shall need a somewhat stronger version of this
inequality.  In particular, note that Lemma \ref{lem:2ndmom} immediately
gives
\begin{equation}
\label{eq:ineq}
\delta(\alpha)<1.78\alpha^2,
\end{equation}
which is stronger than $\delta(\alpha)\le\alpha$ for $\alpha<1/2$.  It is certainly possible
to get improvements on \eqref{eq:ineq} by averaging higher moments of $m/\varphi(m)$,
as was done in \cite{KKP}, which would lead to small improvements on our lower bound
for $C(n)$.

We shall also need some estimates for $\delta(\alpha)$ when $\alpha$
is close to 1, and for this we use an argument of Erd\H os \cite[Theorem 3]{E}.
There he shows, essentially, that $1-\delta(1-\epsilon)\sim2/(e^\gamma|\log\epsilon|)$
as $\epsilon\to0$, where $\gamma$ is Euler's constant.  We will need an estimate with somewhat more precision.

Let
\[
\delta(\alpha,n)=\frac1nD(\alpha,n),\quad
M(x)=\prod_{3\le p\le x}\left(1-\frac1p\right),\quad s_{j,x}=\sum_{4^{j-1}x<p\le 4^jx}\frac1p.
\]
\begin{lemma}
\label{lem:ep}
Uniformly for $2\le x\le\log n$  we have
\[
1-\delta(1-1/x,n)\le M(x)-1/\sqrt{n}
\]
and 
\[
1-\delta(1-1/x,n)\ge M(2x)\left(1-\sum_{j\ge1}\frac{s_{j,2x}^{j+1}}{(j+1)!}\right)
+O\left(\frac1{(\log n)^{\log\log\log n}}\right).
\]
\end{lemma}
\begin{proof}
Note that $1-\delta(1-1/x,n)$ denotes the fraction of numbers $m<2n$ 
 with $\varphi(m)/m > 1-1/x$ (all such $m$ are odd).  In fact, such numbers are not divisible by any prime
$p\le x$, which with Lemma \ref{lem:sieve} gives the upper bound.

For the lower bound we count the numbers $m<2n$ that are not divisible
by any prime $p\le 2x$ and also divisible by at most $j$ distinct primes from
each interval $I_j:=(4^{j-1}\cdot2x,\,4^j\cdot2x]$.   Indeed, if $m$ is such a
number, then
\[
\frac{\varphi(m)}{m}>\prod_{j\ge1}\left(1-\frac1{4^{j-1}\cdot2x}\right)^j
>1-\sum_{j\ge1}\frac j{4^{j-1}\cdot2x}=1-\frac8{9x}.
\]
Let $A_j$ be the set of products of $j+1$ distinct primes from $I_j$.
For $a\in A_j$, the number of odd numbers $m<2n$ with $a\mid m$ and $m$ not
divisible by any prime to $2x$, is by Lemma \ref{lem:sieve}, within
$2^{\pi(2x)-1}$ of $M(2x)n/a$.  Note too that the sum of $1/a$ for $a\in A_j$
is at most $s_{j,2x}^{j+1}/(j+1)!$, by the multinomial theorem.  Let
$\pi(I_j)$ denote the number of primes in $I_j$.  Thus, the number of odd $m<2n$
not divisible by any prime to $2x$ and divisible by some $a\in A_j$ is uniformly
\[
M(2x)ns_{j,2x}^{j+1}/(j+1)!+O\left(2^{\pi(2x)}\binom{\pi(I_j)}{j+1}\right).
\]
The binomial coefficient here is bounded by $O(4^{j(j+1)}x^{j+1})$ using only that 
$\pi(I_j)<4^jx$.  Note that for $j\le\log\log n$ this expression is $O_\epsilon(n^\epsilon)$
for any $\epsilon>0$, as is $2^{\pi(2x)}$, so that the number of integers $m<2n$
not divisible by any prime $p\le2x$ yet divisible by some $a\in A_j$ for
$j\le\log\log n$ is at least
\[
M(2x)n\sum_{j\le\log\log n}\frac{s_{j,2x}^{j+1}}{(j+1)!}+O(n^{1/2}).
\]
Thus,
\begin{align*}
n-D(1-1/x,n)\ge &M(2x)n\left(1-\sum_{j\le\log\log n}\frac{s_{j,2x}^{j+1}}{(j+1)!}\right)\\
&\quad-2n\sum_{j>\log\log n}\frac{s_{j,2x}^{j+1}}{(j+1)!}+O(n^{1/2}).
\end{align*}
Since $s_{j,2x}=O(1/j)$, we have
\[
\sum_{j>\log\log n}\frac{s_{j,2x}^{j+1}}{(j+1)!}=O(1/(\log n)^{\log\log\log n}).
\]
Thus, our count is 
\[
\ge M(2x)n\left(1-\sum_{j\ge1}\frac{s_{j,2x}^{j+1}}{(j+1)!}\right)+O(n/(\log n)^{\log\log\log n}),
\]
which gives our lower bound.
\end{proof}

\begin{corollary}
\label{cor:ep}
Uniformly for $2\le x\le \log n$, we have as $n\to\infty$,
\[
1-\delta(1-1/x,n)\le\frac{2}{e^\gamma\log x}\left(1+\frac1{2(\log x)^2}+o(1)\right).
\]
Further, for $150\le x\le\log n$ and $n$ sufficiently large,
\[
1-\delta(1-1/x,n)\ge\frac2{e^\gamma\log(2x)}\left(1-\frac7{4(\log(2x))^2}\right).
\]
\end{corollary}
\begin{proof}
By Rosser and Schoenfeld \cite[(3.26)]{RS} we have
\[
M(x)<\frac2{e^\gamma\log x}\left(1+\frac1{2(\log x)^2}\right),
\]
so our first assertion follows from the first part of Lemma \ref{lem:ep}.
Further, using \cite[(3.25)]{RS} we have
\[
M(2x)>\frac2{e^\gamma\log(2x)}\left(1-\frac1{2(\log(2x))^2}\right),
\]
and so the second part our our assertion will follow from Lemma \ref{lem:ep}
if we show
\[
\sum_{j\ge1}\frac{s_{j,2x}^{j+1}}{(j+1)!}<\frac {1.4}{(\log(2x))^2}
\]
for all sufficiently large $x$, noting that  $1.4<(2/e^\gamma)1.25$.
Using \cite[(3.17),(3.18)]{RS}, we have
\begin{align*}
s_{j,2x}&<\log\log (2\cdot4^jx)-\log\log(2\cdot4^{j-1}x)+\frac1{(\log(2\cdot4^{j-1}x))^2}\\
&<\frac{\log 4}{\log(2\cdot4^{j-1}x)}+\frac1{(\log(2\cdot4^{j-1}x))^2}
\le \frac{\log 4}{\log(2x)}+\frac1{(\log(2x))^2}=:s.
\end{align*}
Thus, using $x\ge150$,
\[
\sum_{j\ge1}\frac{s_{j,2x}^{j+1}}{(j+1)!}<e^s-1-s<0.55s^2
\]
and $s^2<2.5/(\log(2x))^2$, so our claim follows.
\end{proof}

In addition, we shall use the following numerical bounds.  The first of these
follows from Kobayashi \cite{K}, the last two from Lemma \ref{lem:ep}, and
the others from  Wall \cite{W}.
 \begin{align}
\label{eq:est}
\begin{split}
0.02240&<\delta(0.5)<0.02352,\\
0.1160&<\delta(0.6)<0.1624,\\
0.3556&<\delta(0.7)<0.3794,\\
0.4808&<\delta(0.8)<0.5120,\\
0.5644&<\delta(0.9)<0.6310\\
0.7593&<\delta(0.99)<0.7949\\
0.8380&<\delta(0.999)<0.8539.
\end{split}
\end{align}

\section{The lower bound}

We partition $(0,1]$ into consecutive
intervals 
\[
(\alpha_0,\alpha_1],(\alpha_1,\alpha_2],\dots,(\alpha_{k-1},\alpha_k],
\hbox{ where }
0=\alpha_0<\alpha_1<\dots<\alpha_k=1.
\] 
The parameter $k$ will depend gently on $n$, namely $k=O(\log\log n)$.
The partition of $(0,1]$ will correspond to a partition of $[n]$ into subsets as follows.
For $j=0,1,\dots,k-1$, let
\[
S_j=\big\{m\in\{1,3,\dots,2n-1\}:\alpha_j<\varphi(m)/m\le\alpha_{j+1}\big\}.
\]

In getting a lower bound for $C_0(n)$, we show that  there are many ways to
assign coprime companions for each member $m$ of $\{1,3,\dots,2n-1\}$ that do
not overlap with the choices for other values of $m$.  In particular, we organize
the odd numbers $m<2n$ by increasing size of $\varphi(m)/m$, and so organize
them into the sets $S_1,S_2,\dots$.  In particular, we will choose the parameters
$\alpha_j$ in such a way that there are more ways to assign coprime companions
for $m\in S_j$ than there are members in all of the sets $S_i$ for $i\le j$
combined.

For an odd number $m<2n$ let $F(m,n)$ denote the number of integers
in $[n]$ coprime to $m$.  Suppose $0<\alpha<\beta<1$ and we wish to find
coprime assignments for members of 
\[
S=\{m\hbox{ odd}:m<2n,\,\varphi(m)/m\in(\alpha,\beta]\}=\{m_1,m_2,\dots,m_t\},
\]
where $t=\#S=D(\beta,n)-D(\alpha,n)$.
Let $M=\lceil\alpha n-\sqrt{n}\rceil$, so that for each $m\in S$ we have
$F(m,n)\ge M$, via Corollary \ref{cor:sieve}.  Assume that those odd $m<2n$
with $\phi(m)/m\le \alpha$ already have their coprime assignments.
Then $m_1$ can be assigned to
at least $M-D(\alpha,n)$ numbers in $[n]$, $m_2$ can be assigned to at least
$M-1-D(\alpha,n)$ numbers in $[n]$, etc.  In all, the numbers in $S$ have at
least
\begin{equation}
\label{eq:Scount}
\frac{(M-D(\alpha,n))!}{(M-D(\alpha,n)-\#S)!}
=\frac{(M-D(\alpha,n))!}{(M-D(\beta,n))!}
\end{equation}
coprime assignments that do not interfere with those for $\varphi(m)/m\le\alpha$.
If $0<a<b<1$ and $an,bn$ are integers, then
\[
(an-bn)!=\exp((a-b)n(\log n-1) +(a-b)n\log(a-b)+O(\log n)).
\]
Let $f(x)=x\log x$.
Thus, the expression in \eqref{eq:Scount} is equal to
\[
\exp((\delta(\beta,n)-\delta(\alpha,n))n(\log n-1)+E(\alpha,\beta,n)n+O(\log n)),
\]
where
\[
E(\alpha,\beta,n)=f(\alpha-\delta(\alpha,n))-f(\alpha-\delta(\beta,n)).
\]
We thus will have that
$C_0(n)\ge n!\exp\left(nE+O(k\log n)\right)$, where
\begin{equation}
\label{eq:key}
E=\sum_{1\le i\le k-1}(f(\alpha_i-\delta(\alpha_i,n))-f(\alpha_i-\delta(\alpha_{i+1},n))).
\end{equation}
(We will choose $\alpha_1=1/\log\log n$ and for $n$ sufficiently large,
every odd $m<2n$ will have $\varphi(m)/m>\alpha_1$, so the interval $(0,\alpha_1]$ does
not contribute.)

The sum in \eqref{eq:key} is
almost telescoping.  In particular the density $\delta(\alpha_{i+1},n)$
when $1\le i\le k-2$ appears twice, the two $f$-values being
\[
-f(\alpha_i-\delta(\alpha_{i+1},n))+f(\alpha_{i+1}-\delta(\alpha_{i+1},n)).
\]
We do not have a completely accurate evaluation for $\delta(\alpha_{i+1},n)$ nor
for the limiting value of $\delta(\alpha_{i+1})$, but we do have a fairly narrow
interval where this limit lives.  Note that the expression
\[
-f(\alpha_i-x)+f(\alpha_{i+1}-x)
\]
is decreasing in $x$ when $0<x<\alpha_i$, so if we use an upper bound
for $\delta(\alpha_{i+1},n)$ in \eqref{eq:key}, we will get a lower bound
for the sum.

\subsection{The interval $(0,1/4]$}
Let $j_0$ be the least integer with $2^{j_0}>\log\log n$ and let $\alpha_1=1/2^{j_0}$.
Further, let $\alpha_j=2^{j-1}\alpha_1=1/2^{j_0-j+1}$, for $j\le j_0-1$.  This gives the
first part of our partition of $(0,1]$, namely the sets $(\alpha_{j},\alpha_{j+1}]$ for
$j\le  j_0-1$ give a partition of $(0,1/4]$.  

 Using \eqref{eq:ineq} and the upper bound for $\delta(1/2)$ in \eqref{eq:est}, we
 have
\[
\delta(1/2^i,n)\le\min\{1.78/4^i,\,0.02352\}
\]
for all $i$ and all large $n$.
  We find  the $E$-sum from \eqref{eq:key} for the portion for $(0,1/4]$ is
 $>-0.0538$.  So the contribution
for this part of the count is greater than
\begin{equation}
\label{eq:1part}
\exp\big(D(1/4,n)(\log n-1)-0.0538n\big)
\end{equation}
for all large $n$.

\subsection{The interval $(1/4,0.999]$}
We split the interval $(1/4,0.999]$ at 
\[
0.5,~0.6,~0.7,~0.8,~0.9,~0.99.
\]
Using the upper bounds for our various densities from \eqref{eq:est},
we have the $E$-sum from \eqref{eq:key} is
\begin{align*}
&f(1/4-.02352)-f(1/4-.02352)+f(.5-.02352)-f(0.5-.1624)\\
&+f(.6-.1624)-f(.6-.3794)+f(.7-.3794)-f(.7-.5120)\\
&+f(.8-.5120)-f(.8-.6310)+f(.9-.6310)-f(.9-.7949)\\
&+f(.99-.7949)-f(.99-.8539)
>-0.2873.
\end{align*}
(Note the first two terms are not a typo!)
Thus, the contribution from $(1/4,0.999]$ is greater than

\begin{equation}
\label{eq:2part}
\exp\big((D(0.999,n)-D(1/4,n)(\log n-1)-0.2873n\big)
\end{equation}
for all large $n$.

\subsection{The interval $(0.999,1-1/\log n]$}

Let $j_1$ be the least integer with $10^{j_1}>\log n$.   Let
$\epsilon_i=10^{-i}$.   We 
deal with the intervals 
\[
(1-\epsilon_{i-1},1-\epsilon_i] ~\hbox{ for }~4\le i\le j_1-1.
\]
For our argument to work we will need to show that 
$D(1-\epsilon_i,n)<(1-\epsilon_{i-1})n-\sqrt{n}$, that is,
\begin{equation}
\label{eq:need}
\epsilon_{i-1}n-\sqrt{n}<n-D(1-\epsilon_i,n)\hbox{ for }i\ge4.
\end{equation}
 From Corollary \ref{cor:ep} we have
\[
n-D(1-\epsilon_i,n)>\frac{2n}{e^\gamma\log(2/\epsilon_i)}\left(1-\frac7{4(\log(2/\epsilon_i))^2}\right).
\]
Note that $\log(2/\epsilon_i)=\log2+i\log10$, so that an expression of
magnitude $1/\log(2/\epsilon_i)$ is much larger than $\epsilon_{i-1}$ when $i\ge4$,
so we have \eqref{eq:need}.

We now compute the contribution from the intervals $(1-\epsilon_{i-1},1-\epsilon_i]$
for $i=4,5,\dots,j_1-1$.  This is at least
\[
\exp(D(1-\epsilon_{j_1-1},n)-D(0.999,n)(\log n-1)+En),
\]
where
\[
E=\sum_{4\le i\le j_1-1}f(1-\epsilon_{i-1}-\delta(1-\epsilon_{i-1},n)-f(1-\epsilon_{i-1}-\delta(1-\epsilon_{i},n)).
\]
Using our bound $0.8539$ for $\delta(0.999)$ from \eqref{eq:est} and 
Corollary \ref{cor:ep} for $\delta(1-\epsilon_i)$ for $i\ge4$, we have $E>-0.2814$, so the contribution
for all large $n$ is at least
\[
\exp((D(1-\epsilon_{j_1-1},n)-D(0.999,n))(\log n-1)-0.2814n).
\]
The final interval $(1-\epsilon_{j_1-1},1-1/\log n]$ contributes
\[
\exp(D(1-1/\log n,n)-D(1-\epsilon_{j_1-1},n)(\log n-1)+O(n/\log\log n)),
\]
 so our total
contribution from $(.999,1-1/\log n]$ is at least
\begin{equation}
\label{eq:3part}
\exp((D(1-1/\log n,n)-D(0.999,n)(\log n-1)-0.2815n)
\end{equation}
for all large $n$.

\subsection{The interval $(1-1/\log n,1]$}

We break this interval at $1-1/\sqrt{2n}$.  It is evident that if $m<2n$ is odd and
$\varphi(m)/m>1-1/\sqrt{2n}$, then $m=1$ or $m$ is a prime in the interval $(\sqrt{2n},2n)$.
Thus,
\begin{equation}
\label{eq:top}
D(1-1/\sqrt{2n},n)=n-2n/\log n +O(n/(\log n)^2)
\end{equation}
by the prime number theorem.  Thus, $\delta(1-1/\sqrt{2n},n)<1-1/\log n$ for all large $n$.
  A calculation
shows that the contribution is at least
\[
\exp\left(\Big(D\big(1-\frac1{\sqrt{2n}},n\big)-D\big(1-\frac1{\log n},n\big)\Big)(\log n-1)+E\right),
\]
where $E=O(n\log\log\log n/\log\log n)$,
this term coming from $f(1-1/\log n-\delta(1-1/\log n,n))$.

For the final interval, we have already noted that the numbers in $[n]_\oo$ remaining
are 1 and the primes in $(\sqrt{2n},2n)$.  We follow
the argument in \cite[Proposition 1]{PS}.  Label the primes in $(\sqrt{2n},2n)$
in decreasing order $p_1,p_2,\dots,p_t$, so that, by \eqref{eq:top},
$t= 2n/\log n+O(n/(\log n)^2)$.  Each $p_i$ has $<2n/p_i$ multiples to $2n$,
of which $<n/p_i+1/2$ are odd.  Let $u=\lfloor t/2\rfloor=n/\log n+O(n/(\log n)^2)$, 
so that $p_u\sim n$.   We count assignments
for $p_i$ for $i=t,t-1,\dots,u$ in order.  At each $i$ there are $i+1$ numbers
remaining to be associated with $p_i$ of which at most $n/p_i+1/2$ are
multiples of $p_i$.  So, there are at least $i-\sqrt{n}$ coprime choices for
$p_i$'s assignment.  Multiplying these counts, we have at least
\[
(u-\sqrt{n})^u=\exp(n+O(n/\log n)
\]
choices.  For each of the remaining primes $p_i$ there are $i+1$ numbers left as possible
assignments, with at most one of these divisible by (actually, equal to) $p_i$.
So the contribution of these primes is $(u-1)!=\exp(n+O(n/\log n))$.  The
final number to assign is 1, and it goes freely to the remaining number left.
So for this interval we have at least 
\[
\exp(2n+O( n/{\log n}))
\]
possibilities.  By \eqref{eq:top}
the count can be rewritten as
\[
\exp\left(\big(n-D\big(1-\frac1{\sqrt{2n}},n\big)\big)(\log n-1)+O\big(\frac{n}{\log n}\big)\right)
\]

With the prior calculation, we have at least
\begin{equation}
\label{eq:4part}
\exp((n-D(1-1/\log n,n)(\log n-1)+O(n\log\log\log n/\log n))
\end{equation}
assignments.

To conclude the proof we multiply the expressions in \eqref{eq:1part},
\eqref{eq:2part}, \eqref{eq:3part}, and \eqref{eq:4part}, getting at least
\[
\exp(n(\log n-1)-0.6226n)
\]
coprime matchings from $[n]_\oo$ to $[n]$  for all large $n$.  Since
$e^{0.6226}>1.8637$, this completes the proof of
the lower bound in Theorem \ref{thm:main2}.

\section{The upper bound and a conjecture}
\label{sec:ub}

For each integer $k\ge2$, let $C_k(n)$ denote the number of permutations
$\sigma$ of $[n]$ where $\gcd(j,\sigma(j),k!)=1$ for each $j\in[n]$.  Thus,
$C(n)\le C_k(n)$ for every $k$.  In fact, $C(n)=C_k(n)$ when $k\ge n$,
but we are interested here in the situation when $k$ is fixed and $n$ is
large.  We claim that
for each fixed $k\ge2$ there is a positive constant $c_k$ such that
$C_k(n)=n!/(c_k+o(1))^n$ as $n\to\infty$.  

Here is a possible plan for the proof of this claim.  Let $K$ be the product
of the primes to $k$.
If $ dd' \mid K$, then one can count the number of $m\in[n]$ with $\gcd(m,K)=d$
that get mapped to an $m'$ with $\gcd(m',K)=d'$.  The product of all of the positive
counts is $n^{O(1)}$, so basically, up to a factor of this shape, the
number of permutations is given by those with one optimal suite of counts.

Let $I_d$ be the set of $m$ with $\gcd(m,K)=d$ and let $\beta(d,d')$ be the
proportion of members of $I_d$ that get sent to $I_{d'}$ by a given
permutation.  Then for a fixed $d$, the numbers $\beta(d,d')$ have sum 1 for 
$d' \mid K/d$, and sum 1 for a fixed $d'$ and $d\mid K/d' $.  One can 
start with some suite of proportions $\beta(d,d')$ that are
``legal" and consider permutations which approximate these proportions,
and see the count as some complicated, but continuous function of the
variables $\beta(d,d')$.  So, there is an optimal suite of proportions,
via calculus, and this gives rise to $c_k$.

Assume that $c_k$ exists.
Note that the sequence $(c_k)$ is monotone nondecreasing and
that if $p<p'$ are consecutive primes, then $c_k=c_p$ for $p\le k<p'$.
It follows from our lower bound for $C(n)$ that the numbers $c_k$
are bounded above.  Let $c_0=\lim_{k\to\infty}c_k$.
\begin{conjecture}
\label{conj:main}
We have $C(n)=n!/(c_0+o(1))^n$ as $n\to\infty$.
\end{conjecture}

We now prove for $k=2,3,5$ that $c_k$ exists and we compute it.
Our value for $c_5$ gives our upper bound theorem for $C(n)$.

The results in Section \ref{sec:pre} largely carry over in the case $k=2$.
Indeed, note that $C_0(n)\le n!$ and $C_1(n)\le (n+1)!$, so that $C(2n)\le n!^2$
and $C(2n+1)\le (n+1)!^2$.  From this we immediately get that
 $C_2(n)\le n!/(2+o(1))^n$ as $n\to\infty$.  In fact, from the proof of Lemma
 \ref{lem:odd} we have $C_2(2n)=n!^2$ and $C_2(2n+1)=(n+1)!^2$,
 so that $c_2=2$.

For $k=3$, we first deal with $6n$ and count one-to-one functions $\sigma$ from $\{1,2,\dots,3n\}$ to $\{1,3,\dots,6n-1\}$ that map multiples of 3 to non-multiples of 3.
There are precisely $(2n)!^2/n!$ of them, so $C_3(6n)=((2n)!^2/n!)^2$.  Similarly
we get $C_3(6n+3)=((2n+1)!^2/(n+1)!)^2$, so these two formulas lead to
$C_3(n)=n!/(3/2^{1/3}+o(1))$ as $n\to\infty$ with $3\mid n$.  To get to other
cases, note that $C_3(n)\le C_3(n+2)$ for all $n$, so we can sandwich $n$
between 2 consecutive multiples of 3 and absorb the error in the ``$o(1)$".
We thus have
\[
c_3=2^{-1/3}3=2.381101\dots.
\]

The case $k=5$ is considerably harder.  We only treat multiples
of 30, the case 15 (mod 30) is similar, and since $C_k(n)\le C_k(n+2)$, we can extend
to all $n$ readily.  The problem is reduced to counting matchings
from $[15n]$ to $[15n]_\oo$ where corresponding terms have gcd coprime
to 15.  We split $[15n]$ into the $n$ multiples of 15, the $2n$ numbers that are
divisible by 5 but not 3, the $4n$ numbers divisible by 3 but not 5, and the
$8n$ numbers coprime to 15.  We have the corresponding decomposition
for $\{1,3,\dots,30n-1\}$.  The first group consisting of the multiples of 15
must be mapped to the numbers coprime to 15, and this can be done in
\[
\frac{(8n)!}{(7n)!}
\]
ways.  The next case we consider is the $2n$ multiples of 5 but not 3.
They must be mapped to the numbers coprime to 5, where some of them
are mapped to numbers coprime to 15 and the rest of them are mapped
to numbers divisible by 3 but not 5.  A calculation shows that the most numerous
case is when it is half and half, also considering the next step which is to
place the multiples of 3 but not 5.  So the total will be within a factor $2n$
of this most numerous case, which has
\[
\binom{2n}{n}\frac{(7n)!}{(6n)!}\frac{(4n)!}{(3n)!}
\]
matchings.  For the multiples of 3 but not 5, these are mapped into the union
of the remaining $6n$ numbers coprime to 15 and the $2n$ numbers divisible
by 5 but not 3, for a total of
\[
\frac{(8n)!}{(4n)!}
\]
matchings.  The remaining $8n$ numbers are all coprime to 15 and can
be mapped to the remaining $8n$ numbers in every possible way, giving
$(8n)!$ matchings.  In all we thus have
\[
\frac{(8n)!}{(7n)!}\binom{2n}{n}\frac{(7n)!}{(6n)!}\frac{(4n)!}{(3n)!}\frac{(8n)!}{(4n)!}(8n)!n^{O(1)}
=\frac{(8n)!^3(2n)!}{(6n)!(3n)!n!^2}n^{O(1)}
\]
matchings.  The log of this expression is within $O(\log n)$ of
\[
=\exp(15n(\log n-1)+(24\log 8+2\log2-6\log6-3\log3)n).
\]
Our count is then squared and $(30n)!$ is factored out, giving
\[
(30n)!\exp((136\log2-18\log3-30\log30)n+O(\log n)).
\]
This then gives that
\[
C_5(n)=(n!/c_5^n)n^{O(1)},
\]
where 
\[
c_5=\exp\big(-\frac{53}{15}\log2+\frac85\log3+\log5\big)
=2^{-53/15}3^{8/5}5=2.504521\dots.
\]

\subsection{A possible value for $c_0$}
\label{sec:ubmcnew}
Nathan McNew has suggested the following argument.  First, for a prime $p$, let
$N_p(n)$ be the number of permutations $\sigma$ of $[n]$ with each $\gcd(j,\sigma(j),p)=1$.
So the constraint is that the multiples of $p$ get mapped to the non-multiples of $p$, and
so we have
\[
N_p(n)=\frac{(\lfloor(1-1/p)n\rfloor!)^2}{\lfloor(1-2/p)n\rfloor!}n^{O(1)}.
\]
Then, up to a factor $n^{O(1)}$, we have
\[
\frac{n!}{N_p(n)}=\left(\frac{p(p-2)^{1-2/p}}{(p-1)^{2(1-1/p)}}\right)^n,
\]
which suggests by independence that
\[
c_k=\prod_{p\le k}\frac{p(p-2)^{1-2/p}p}{(p-1)^{2(1-1/p)}}.
\]
(We interpret the factor at $p=2$ as 2.)  This expression agrees with our
computation of $c_k$ for $k$ up to 5.  And it suggests that $c_0$ is the infinite
product over all primes $p$, so that $c_0=2.65044\dots$.

\section{Computing $C(n)$}
\label{sec:comp}

In this section we discuss the numerical computation of $C(n)$ for modest
values of $n$.  In \cite{O} it is remarked that $C(n)$ has been computed to
$n=30$ by Seiichi Manyama, and extended to $n=50$ by Stephen Locke,
see https://oeis.org/A005326/b005326.txt.  We have verified these values
using the methods of this section and Mathematica.

\begin{footnotesize}
\begin{table}[]
\caption{Values of $C_0(n)=\sqrt{C(2n)}$ and $r_{2n}$.}
\label{Ta:C0}
\begin{tabular}{|rrrr|} \hline
$n$&$C_0(n)$ & $r_{2n}$&\\ \hline
1 & 1 & 1.4142&\\
2 & 2 & 1.5651&\\
3 & 4 & 1.8860&\\
4 & 18 & 1.8276&\\
5 & 60 & 1.9969&\\
6 & 252 & 2.1044&\\
7 & 1{,}860 & 2.0625&\\
8 & 9{,}552 & 2.1629&\\
9 & 59{,}616 & 2.2260&\\
10 & 565{,}920 & 2.2082&\\
11 & 4{,}051{,}872 & 2.2707&\\
12 & 33{,}805{,}440 & 2.3118&\\
13 & 465{,}239{,}808 & 2.2727&\\
14 & 4{,}294{,}865{,}664 & 2.3171&\\
15 & 35{,}413{,}136{,}640 & 2.3850&\\
16 & 768{,}372{,}168{,}960 & 2.3122&\\
17 & 8{,}757{,}710{,}173{,}440 & 2.3451&\\
18 & 79{,}772{,}814{,}777{,}600 & 2.4122&\\
19 & 1{,}986{,}906{,}367{,}584{,}000 & 2.3531&\\
20 &  22{,}082{,}635{,}812{,}268{,}800 & 2.4029&\\
21 & 280{,}886{,}415{,}019{,}776{,}000 & 2.4374&\\
22 & 7{,}683{,}780{,}010{,}315{,}046{,}400 & 2.3905&\\
23 & 102{,}400{,}084{,}005{,}498{,}547{,}200 & 2.4278&\\
24 & 1{,}774{,}705{,}488{,}555{,}494{,}476{,}800 & 2.4401&\\
25 & 40{,}301{,}474{,}964{,}335{,}327{,}232{,}000 & 2.4291&\\
\hline
\end{tabular}
\end{table}
\end{footnotesize}

As is easy to see, the permanent of the incidence matrix of a bipartite
graph of two $n$-sets gives the number of perfect matchings contained
in the graph.  Let ${\bf B}(n)$ be the $n\times n$ ``coprime matrix", where
${\bf B}(n)_{i,j}=1$ when $i,j$ are coprime and 0 otherwise.
So, in particular, and as noted by Jackson \cite{J},
\begin{equation}
\label{eq:perm}
C(n)=\perm({\bf B}(n)).
\end{equation}
However, it is not so simple to compute a large permanent, though we do have
some algorithms that are better than brute force, for example \cite{R} and \cite{B}.

Recall from Lemma \ref{lem:odd} that $C(2n)=C_0(n)^2$, where $C_0(n)$ is
the number of coprime matchings between $[n]$ and $[n]_\oo$.  Thus,
$C(2n)$ can be obtained from an $n\times n$ permanent, which is
considerably easier than the more naive $2n\times2n$ permanent required 
when applying \eqref{eq:perm} to $C(2n)$.

There is a similar reduction for computing $C(2n+1)$.  For each $a\in[n+1]_\oo$,
let $C_{(a)}(n)$ denote the number of coprime matchings between $[n]$ and
$[n+1]_\oo\setminus\{a\}$.  Then $C_1(n)=\sum_{a\in[n+1]_\oo}C_{(a)}(n)$ and
\[
C(2n+1)=\sum_{a\in[n+1]_\oo}\sum_{\substack{b\in[n+1]_\oo\\\gcd(a,b)=1}}
C_{(a)}(n)C_{(b)}(n).
\]
Thus, $C(2n+1)$ can be easily computed from $n+1$ permanents of size $n\times n$.

Let $r_n=(n!/C(n))^{1/n}$, so that $C(n)=n!/r_n^n$.  We have shown that for all
large $n$ we have $2.5<r_n<3.73$.  In the following tables we have computed
the actual values of $r_n$ for $n\le 50$ rounded to 4 decimal places.

It is easy to see that $C_0(n)$ is the number of partitions of $[2n]$ into coprime
unordered pairs.  This has its own OEIS page: A009679, and has been enumerated
there up to $n=30$.

\begin{footnotesize}
\begin{table}[]
\caption{Values of $C(n)$ for $n$ odd and $r_{n}$.}
\label{Ta:C1}
\begin{tabular}{|rrrr|} \hline
$n$&$C(n)$ & $r_{n}$&\\ \hline
1 & 1 & 1&\\
3 & 3 & 1.2599&\\
5 & 28 & 1.3378&\\
7 & 256 & 1.5307&\\
9 & 3{,}600 &1.6696&\\
11 & 129{,}774 & 1.6834&\\
13 & 3{,}521{,}232 & 1.7776&\\
15 & 60{,}891{,}840 & 1.9444&\\
17 & 8{,}048{,}712{,}960 & 1.8761&\\
19 & 425{,}476{,}094{,}976 & 1.9372&\\
21 & 12{,}474{,}417{,}291{,}264 & 2.0648&\\
23 &2{,}778{,}580{,}249{,}611{,}264 & 2.0090&\\
25 & 172{,}593{,}628{,}397{,}420{,}544 & 2.0804&\\
27 & 17{,}730{,}530{,}614{,}153{,}986{,}048 & 2.1159&\\
29 & 4{,}988{,}322{,}633{,}552{,}214{,}818{,}816 & 2.0841&\\
31 & 427{,}259{,}978{,}841{,}815{,}654{,}400{,}000 & 2.1466&\\
33 & 57{,}266{,}563{,}000{,}754{,}880{,}493{,}977{,}600 & 2.1818&\\
35 & 14{,}786{,}097{,}120{,}330{,}296{,}843{,}693{,}260{,}800 & 2.1798&\\
37 & 3{,}004{,}050{,}753{,}199{,}657{,}126{,}879{,}764{,}480{,}000 & 2.1988&\\
39 & 536{,}232{,}134{,}065{,}318{,}935{,}894{,}365{,}552{,}640{,}000 & 2.2295&\\
41 & 274{,}431{,}790{,}155{,}416{,}580{,}402{,}144{,}584{,}785{,}920{,}000 & 2.2058&\\
43 & 51{,}681{,}608{,}012{,}142{,}138{,}983{,}265{,}921{,}023{,}262{,}720{,}000 & 2.2409&\\
45& 7{,}417{,}723{,}304{,}411{,}612{,}192{,}092{,}096{,}851{,}178{,}291{,}200{,}000 & 2.2918&\\
47 & 7{,}896{,}338{,}788{,}322{,}918{,}879{,}731{,}318{,}625{,}512{,}774{,}041{,}600{,}000 & 2.2459&\\
49 & 1{,}989{,}208{,}671{,}980{,}285{,}257{,}956{,}064{,}090{,}726{,}080{,}876{,}380{,}160{,}000 & 2.2743&\\
\hline
\end{tabular}
\end{table}
\end{footnotesize}

\section{Anti-coprime permutations}
\label{sec:anti}

One might also wish to consider permutations $\sigma$ of $[n]$ where each $\gcd(j,\sigma(j))>1$.
Of course, none exist, since $1\in[n]$.  Instead we can count the number $A(n)$ where
$\gcd(j,\sigma(j))>1$ for $2\le j\le n$.  This seems like an interesting problem.  We can
prove the following lower bound.

\begin{proposition}
\label{prop:anti}
As $n\to\infty$, we have
\begin{equation}
\label{eq:anti}
A(n) \ge n!/\exp((e^{-\gamma}+o(1))n\log\log n).
\end{equation}
\end{proposition}

We sketch the proof.  Let $\epsilon_n=1/\sqrt{\log\log n}$ and let $g(x)=\prod_{p<x}(1-1/p)$,
so that $g(x)$ is similar to the function $M(x)$ we considered earlier.  For each prime $p<n^{\epsilon_n}$ consider 
the set $L_n(p)$ of integers $m\le n$ with least prime factor $p$, and let
\[
\lambda(p,n)=\frac1n\#L_n(p).
\]
Note that 
\[
\bigcup_{p<n^{\epsilon_n}}L_n(p)
\]
is the set of integers with least prime factor $<n^{\epsilon_n}$, so the number of
integers $m\le n$ not in this union is $O(n/(\epsilon_n\log n))$.  In particular,
\begin{equation}
\label{eq:lambda}
\sum_{p<n^{\epsilon_n}}\lambda(p,n)=1+O(1/(\epsilon_n\log n)).
\end{equation}

Each of the $(\#L_n(p))!$ permutations of $L_n(p)$ is anti-coprime, and gluing these together
for $p<n^{\epsilon_n}$ and having the remaining
elements of $[n]$ as fixed points, gives an anti-coprime permutation of $[n]$.  
So we have
\[
A(n) \ge \prod_{p<n^{\epsilon_n}}(\#L_n(p))!.
\]
Thus, by the inequality $k!>(k/e)^k$,
\begin{align*}
A(n)&\ge\exp\left(\sum_{p<n^{\epsilon_n}}\lambda(p,n)n(\log n+\log(\lambda(p,n))-1)\right)\\
&=\exp\left(\sum_{p<n^{\epsilon_n}}\lambda(p,n)n(\log n-1)+nE\right)
\end{align*}
where
\[
E =\sum_{p<n^{\epsilon_n}}\lambda(p,n)\log(\lambda(p,n)).
\]

Note that by \eqref{eq:lambda}
\[
\sum_{p<n^{\epsilon_n}}\lambda(p,n)n(\log n-1)=n\log n+O(n/\epsilon_n).
\]

 To deal with $E$, we have that for $p<n^{1/\epsilon_n}$,
\begin{equation*}
\lambda(p,n)=\frac{ g(p)}{p}(1+O(e^{-1/\epsilon_n }))
\end{equation*}
uniformly for large $n$.  Indeed each $m\in L_n(p)$ is of the form $pk$ where $k\le n/p$
is an integer not divisible by any prime $q< p$.  Such integers $k$ are easily counted
by the fundamental lemma of either Brun's or Selberg's sieve, which gives the above estimate.  

We have $g(p)$ of magnitude $1/\log p$,  in fact $g(p)=1/(e^\gamma\log p)(1+O(1/\log p))$.
Thus, we have
\begin{align*}
E&=\sum_{p<n^{\epsilon_n}}\frac{g(p)}{p}(\log(g(p)/p)(1+O(e^{-1/\epsilon_n }))\\
&\kern-3pt=\sum_{p<n^{\epsilon_n}}\frac1{e^\gamma p\log p}(-\log p-\log\log p-\gamma)(1+O(1/\log p)
+O(e^{-1/\epsilon_n }))\\
&\kern-3pt=-\sum_{p<n^{\epsilon_n}}\frac1{e^\gamma p}(1+O(e^{-1/\epsilon_n })+O(1).
\end{align*}
It remains to note that 
\[
\sum_{p<n^{\epsilon_n}}\frac1p=\log\log n-\log\epsilon_n +O(1).
\]
Thus, we have Proposition \ref{prop:anti}.

We conjecture that $A(n)=n!/\exp((e^{-\gamma}+o(1))n\log\log n)$ as $n\to\infty$,
that is, Proposition \ref{prop:anti} is best possible.  Though it is difficult to ``see" $\log\log n$
tending to infinity, we have some scant evidence in
Table \ref{Ta:A}.  Let $u_n=(n!/A(n))^{1/n}$, so the conjecture is that
$u_n\sim e^{-\gamma}\log\log n$.

The computation of $A(n)$ is helped by the realization that all of the permutations counted
have 1 and the primes in $(n/2,n]$ as fixed points, so one can deal with a somewhat
smaller adjacency matrix 
than $n\times n$.  In particular, if $n$ is prime, then $A(n)=A(n-1)$, so in Table \ref{Ta:A}
we only consider $n$ composite (the cases $n=1,2$ being trivial).
 In addition, for a prime $p\in(n/3,n/2]$ either $p$ is a fixed point
or $(p,2p)$ is a 2-cycle, which gives another reduction.  

\begin{footnotesize}
\begin{table}[]
\caption{Values of $A(n)$ for $n$ composite and $u_{n}$.}
\label{Ta:A}
\begin{tabular}{|rrrr|} \hline
$n$&$A(n)$ & $u_{n}$&\\ \hline
4&2&1.8612&\\
6&8&2.1170&\\
8&30&2.4607&\\
9&72&2.5786&\\
10&408&2.4826&\\
12&4{,}104&2.6440&\\
14&29{,}640&2.8976&\\
15&208{,}704&2.8388&\\
16&1{,}437{,}312&2.8034&\\
18& 22{,}653{,}504&2.9479&\\
20&318{,}695{,}040&3.1199&\\
21&2{,}686{,}493{,}376&3.0866&\\
22&27{,}628{,}410{,}816&3.0356&\\
24&575{,}372{,}874{,}240&3.1722&\\
25&1{,}775{,}480{,}841{,}216&3.2935&\\
26&21{,}115{,}550{,}048{,}256&3.2420&\\
27&132{,}879{,}856{,}582{,}656&3.2758&\\
28&2{,}321{,}256{,}928{,}702{,}464&3.1932&\\
30&83{,}095{,}013{,}944{,}442{,}880&3.2870&\\
\hline
\end{tabular}
\end{table}
\end{footnotesize}

   


 
   
\subsection{Other types of permutations}

One might consider other arithmetic constraints on permutations.  For example,
what can be said about the number of permutations $\sigma$ of $[n]$ where   for each
$j\in[n]$, either $j\mid\sigma(j)$ or $\sigma(j)\mid j$?  Or, the number where each
lcm$[j,\sigma(j)]\le n$?  Problems such as the longest possible cycle in such permutations, 
the minimum number of disjoint cycles, etc.\ were studied in
\cite{P2}, \cite{S}, \cite{M} and elsewhere.  The enumeration problems have not been
well-studied, though the first one has an OEIS page: A320843.
      
\section*{Acknowledgments}
I thank Sergi Elizalde for informing me of \cite{J} and \cite{O} and I am
grateful to Nathan McNew for suggesting the argument in Section \ref{sec:ubmcnew}.

\end{document}